\documentclass[10pt]{amsart}
\input{diagrams}

\newtheorem{proposition}{Proposition}[section]
\newtheorem{lemma}[proposition]{Lemma}
\newtheorem{corollary}[proposition]{Corollary}
\newtheorem{theorem}[proposition]{Theorem}

\theoremstyle{definition}

\theoremstyle{remark}
\newtheorem{remark}[proposition]{Remark}

\def\leftact{\hbox{$\rightharpoonup$}}
\def\rightact{\hbox{$\leftharpoonup$}}

\def\equal#1{\smash{\mathop{=}\limits^{#1}}}

\newcommand{\thlabel}[1]{\label{th:#1}}
\newcommand{\thref}[1]{Theorem~\ref{th:#1}}
\newcommand{\selabel}[1]{\label{se:#1}}
\newcommand{\seref}[1]{Section~\ref{se:#1}}
\newcommand{\lelabel}[1]{\label{le:#1}}
\newcommand{\leref}[1]{Lemma~\ref{le:#1}}
\newcommand{\prlabel}[1]{\label{pr:#1}}
\newcommand{\prref}[1]{Proposition~\ref{pr:#1}}
\newcommand{\colabel}[1]{\label{co:#1}}
\newcommand{\coref}[1]{Corollary~\ref{co:#1}}
\newcommand{\relabel}[1]{\label{re:#1}}
\newcommand{\reref}[1]{Remark~\ref{re:#1}}

\newcommand{\eqlabel}[1]{\label{eq:#1}}
\newcommand{\equref}[1]{(\ref{eq:#1})}

\newcommand{\Hom}{{\rm Hom}}

\newcommand{\End}{{\rm End}}

\newcommand{\Ker}{{\rm Ker}\,}

\newcommand{\tr}{{\rm tr}\,}

\def\lan{\langle}
\def\ran{\rangle}
\def\ot{\otimes}

\newcommand{\Cc}{\mathcal{C}}

\newcommand{\Mm}{\mathcal{M}}

\def\text#1{{\rm {\rm #1}}}

\def\ol{\overline}
\def\ul{\underline}

\begin{document}
\title[Fully bounded noetherian rings]{Fully bounded noetherian rings and
Frobenius extensions}
\author{S. Caenepeel}
\address{Faculty of Engineering,
Vrije Universiteit Brussel, VUB, B-1050 Brussels, Belgium}
\email{scaenepe@vub.ac.be}
\urladdr{http://homepages.vub.ac.be/\~{}scaenepe/}
\author{T. Gu\'ed\'enon}
\address{Faculty of Engineering,
Vrije Universiteit Brussel, VUB, B-1050 Brussels, Belgium}
\email{tguedeno@vub.ac.be, guedenon@caramail.com}
\thanks{Research supported by the project G.0278.01 ``Construction
and applications of non-commutative geometry: from algebra to physics"
from FWO Vlaanderen}
\begin{abstract}
Let $i:\ A\to R$ be a ring morphism, and $\chi:\ R\to A$ a right $R$-linear map
with $\chi(\chi(r)s)=\chi(rs)$ and $\chi(1_R)=1_A$. If $R$ is a Frobenius $A$-ring, then
we can define a trace map $\tr:\ A\to A^R$. If there exists an element of trace 1
in $A$, then $A$ is right FBN if and only if $A^R$ is right FBN and $A$ is
right noetherian. The result can be generalized to the case where $R$ is an
$I$-Frobenius $A$-ring. We recover results of Garc\'{\i}a and del R\'{\i}o and by
D\v{a}sc\v{a}lescu, Kelarev and Torrecillas on actions of group and Hopf algebras
on FBN rings as special cases. We also obtain applications to extensions
of Frobenius algebras, and to Frobenius corings with a grouplike element.
\end{abstract}

\maketitle
\section*{Introduction}
A ring $A$ is called right bounded if every essential right ideal contains a non-zero two-sided ideal.
$A$ is right fully bounded noetherian or right FBN if $A$ is noetherian, and
$A/P$ is right bounded for every two-sided prime ideal $P$ of $A$.
Obviously commutative noetherian rings are right FBN; more generally,
noetherian PI-rings and artinian rings are FBN. A series of conjectures in classical
ring theory can be proved in the case of rings with the FBN property, we refer
to the introduction of \cite{Sorin} for a brief survey.\\
Assume that a finite group $G$ acts on $A$. Garc\'{\i}a and Del R\'{\i}o \cite{Garcia}
investigated the relationship between the FBN property for $A$ and its subring
of invariants $A^G$. The main result is that, in case $A$ is right noetherian,
the right FBN property for $A$ is equivalent to the right FBN property for $A^G$,
if there exists an element in $A$ having trace $1$. A similar statement was proved
in \cite{Nasta} for rings graded by a finite group $G$.
These results can be generalized to Hopf algebra actions (see
\cite{Sorin,Guedenon}). \\
We have observed that the methods introduced in \cite{Garcia} can be applied in an apparently completely
different situation. Let $S$ be a Frobenius algebra (with Frobenius system
$(e=e^1\ot e^2,\ol{\nu})$) and $j:\ S\to A$ an algebra map, with $A$ a right
noetherian ring. If there exists $a\in A$ such that $j(e^1)aj(e^2)=1$, then
$A$ is right FBN if and only if $C_S(A)$ is right FBN.\\
In this note, we propose a unified approach to these results, based on the concept
of an $A$-ring with a grouplike character, as introduced in \cite{CVW}. Basically, this
consists of a ring morphism $i:\ A\to R$, together with a right $A$-linear map
$\chi:\ R\to A$ such that the formula $a\rightact r=\chi(ar)$ makes $A$ into a
right $R$-module. The subring of invariants is defined as $B=\{b\in A~|~b\chi(r)=\chi(br)\}$.
The main result is basically  the following: if $R$ is a Frobenius $A$-ring, and $A$ is projective
as a right $R$-module, then $A$ is right FBN if and only if $B$ is right FBN and
$A$ is right noetherian. The methods of proof are essentially the same as in
\cite{Garcia}. If $R$ is a Frobenius $A$-ring, then we can define a trace map
$\tr:\ A\to B$, and $A$ is projective (and a fortiori quasi-projective) as a right
$R$-module if and only if there exists an element of trace 1. The condition that
$R$ is Frobenius can be relaxed in the sense that it suffices that $R$ is
Frobenius of the second kind, with respect to a strict Morita context 
$(A,A,I,J,f,g)$. Then the trace map is a map $\tr:\ J\to B$.\\
The above mentioned results on group and Hopf algebra actions and extensions
of Frobenius algebras can be obtained as special cases. We also present an
application to Frobenius corings with a grouplike element.

\section{Rings with a grouplike character}\selabel{1}
Let $A$ be an associative ring with unit. The category of $A$-bimodules
${}_A\Mm_A$ is a monoidal category, and we can consider algebras in
${}_A\Mm_A$. Such an algebra $R$ is a ring $R$ together with a ring
morphism $i:\ A\to R$. The bimodule structure on $A$ is then given by
$arb=i(a)ri(b)$, for all $a,b\in A$ and $r\in R$.
A {\sl right grouplike character} on $R$ is a right $A$-linear map
$\chi:\ R\to A$ such that
\begin{equation}\eqlabel{1.1.0}
\chi(\chi(r)s)=\chi(rs)~~{\rm and}~~\chi(1_R)=1_A,
\end{equation}
for all $r,s\in R$. We then say that $(R,i,\chi)$ is an $A$-ring with
a right grouplike character.
Right grouplike characters were introduced in \cite{CVW}.
The terminology is motivated by the fact that the dual of a coring with
a grouplike element is a ring with a grouplike character (see \seref{7}).
For all $a\in A$, we have that
$$\chi(i(a))=\chi(1_R\cdot a)=\chi(1_R)a=1_Aa=a,$$
so $\chi\circ i={\rm Id}_A$, and $i$ is injective, $\chi$ is surjective. Sometimes we will regard $i$ as
an inclusion.
$A$ is a right $R$-module, with right $R$-action
\begin{equation}\eqlabel{1.1.1}
a\rightact r=\chi(ar).
\end{equation}
$A$ is a cyclic right $R$-module, since
$$a=\chi(i(a))=\chi(1_Ai(a))=1_A\leftact i(a),$$
for all $a\in A$. For $M\in \Mm_R$, the submodule of invariants is defined as
$$M^R=\{m\in M~|~mr=m\chi(r),~{\rm for~all~}r\in R\}.$$
Let
$$B=A^R=\{b\in A~|~b\chi(r)=\chi(br),~{\rm for~all~}r\in R\}.$$
Then $B$ is a subring of $A$, $M^R$ is a right $B$-module, and we have the
invariants functor $(-)^{R}:\ \Mm_R\to \Mm_B$. We will now present some
elementary properties of
$$Q=R^R=\{q\in R~|~qr=q\chi(r),~{\rm for~all~}r\in R\}.$$

\begin{lemma}\lelabel{1.1}
Let $(R,i,\chi)$ be an $A$-ring with a right grouplike character.
\begin{enumerate}
\item $Q$ is a $(R,B)$-subbimodule of $R$;
\item $\chi$ restricts to a $B$-bimodule map $\chi:\ Q\to B$;
\item if $1_R\in Q$, then $i$ is an isomorphism of rings, with inverse $\chi$.
\end{enumerate}
\end{lemma}

\begin{proof}
1) We refer to \cite[Prop. 2.2]{CVW}.\\
2) For all $q\in Q$ and $r\in R$, we have
$$\chi(q)\chi(r)=\chi(q\chi(r))=\chi(qr)=\chi(\chi(q)r),$$
hence $\chi(q)\in B$. $\chi$ is right $A$-linear, so its restriction to $Q$
is right $B$-linear. For all $q\in Q\subset R$ and $b\in B$, we have,
by the definition of $A^R=B$ that $b\chi(q)=\chi(bq)$, so $\chi$ is also
left $B$-linear.\\
3) If $1_R\in Q$, then we have for all $r\in R$ that
$$r=1_Rr=1_R\chi(r)=1_Ri(\chi(r))=i(\chi(r)).$$
It follows that $i$ is a left inverse of $\chi$. We have seen above that $i$
is always a right inverse of $\chi$, so it follows that $i$ is an isomorphism.
\end{proof}

If $M\in \Mm_R$, then $\Hom_R(A,M)\in \Mm_B$, with right $B$-action
$(fb)(a)=f(ba)$, for all $b\in B$, $f\in \Hom_R(A,M)$ and $a\in A$.\\
$\End_R(A)$ is a $B$-bimodule, with left $B$-action $(bf)(a)=bf(a)$,
for all $b\in B$, $f\in \End_R(A)$ and $a\in A$.

\begin{lemma}\lelabel{1.2}
Let $(R,i,\chi)$ be an $A$-ring with a right grouplike character,
and $M$ a right $R$-module.
\begin{enumerate}
\item $\Hom_R(A,M)\cong M^R$ as right $B$-modules;
\item $\End_R(A)\cong B$ as $B$-bimodules and as rings.
\end{enumerate}
\end{lemma}

\begin{proof}
1) For $f\in \Hom_R(A,M)$ and $r\in R$, we have
$$f(1_A)r=f(1_A\rightact r)=f(\chi(r))=f(1_A)\chi(r),$$
so $f(1_A)\in M^R$, and we have a well-defined map
$$\phi:\ \Hom_R(A,M)\to M^R,~~\phi(f)=f(1_A).$$
$\phi$ is right $B$-linear since
$$\phi(fb)=(fb)(1_A)=f(b1_A)=f(1_Ab)=f(1_A)b=\phi(f)b.$$
The inverse of $\phi$ is given by the formula
$$\phi^{-1}(m)(a)=ma,$$
for all $m\in M^R$ and $a\in A$.\\
2) If $M=A$, then $\phi$ is also left $B$-linear since
$$\phi(bf)=(bf)(1_A)=bf(1_A)=\phi(f)b.$$
\end{proof}

\section{Quasi-projective modules}\selabel{2}
A right $R$-module $M$ is called {\sl quasi-projective} if the canonical
map $\Hom_R(M,M)\to \Hom_R(M,M/N)$ is surjective, for every $R$-submodule
$N$ of $M$. This means that 
every right $R$-linear map $f:\ M\to M/N$
factorizes through the canonical projection $p:\ M\to M/N$, that is,
there exists a right $R$-linear map $g:\ M\to M$ such that
$f=p\circ g$.

\begin{proposition}\prlabel{2.1}
Let $(R,i,\chi)$ be an $A$-ring with a right grouplike character. The following
assertions are equivalent.
\begin{enumerate}
\item $A$ is quasi-projective as a right $R$-module;
\item for every right $R$-submodule $I$ of $A$, and every $a+I\in (A/I)^R$,
there exists $b\in B$ such that $b-a\in I$;
\item for every right $R$-submodule $I$ of $A$, $(A/I)^R\cong (B+I)/I$.
\end{enumerate}
\end{proposition}

\begin{proof}
$\ul{1)\Rightarrow 2)}$. Observe that
\begin{equation}\eqlabel{2.1.1}
(A/I)^R=\{a+I\in A/I~|~a\chi(r)-\chi(ar)\in I,~{\rm for~all~}r\in R\}.
\end{equation}
For $a+I\in (A/I)^R$, we have a well-defined right $A$-linear map
$$f:\ A\to A/I,~~f(a')=aa'+I.$$
$f$ is right $A$-linear since
\begin{eqnarray*}
&&\hspace*{-2cm} f(a'\rightact r)=a(a'\rightact r)+I=a\chi(a'r)+I\\
&=&\chi(aa'r)+I=((aa')\rightact r)+I=f(a')\rightact r.
\end{eqnarray*}
Let $p:\ A\to A/I$ be the canonical projection. Since $A$ is quasi-projective,
there exists $g\in \Hom_R(A,A)$ such that $p\circ g= f$, that is
$aa'+I=g(a')+I$ and, in particular, $a+I=g(1_A)+I$, or $g(1_A)-a\in I$.
Let us show that $b=g(1_A)\in B$. Indeed, for all $r\in R$, we have
\begin{eqnarray*}
&&\hspace*{-2cm}
\chi(br)-b\chi(r)=\chi(g(1_A)r)-g(1_A)\chi(r)
= (g(1_A)\rightact r)-(g(1_A)\rightact (i\circ \chi)(r))\\
&=& g(1_A\rightact r)-g((\chi\circ i\circ \chi)(r))
= g(\chi(r))-g(\chi(r))=0.
\end{eqnarray*}
$\ul{2)\Rightarrow 3)}$. The map $B\to (A/I)^R$, $b\mapsto b+I$ induces a
monomorphism $(B+I)/I\to (A/I)^R$. Condition 2) means precisely that this
map is surjective.\\
$\ul{3)\Rightarrow 1)}$. Take a right $R$-linear map $f:\ A\to A/I$, with
$I$ a right $R$-submodule of $A$. Then
$$\chi(f(1_A)r)=f(1_A)\rightact r=f(1_A\rightact r)=
f(\chi(1_Ar))=f(\chi(r))=f(1_A)\chi(r),$$
so $f(1_A)\in (A/I)^R\cong (B+I)/I$. Take $b\in B$ such that
$f(1_A)=b+I$, and consider the map $g:\ A\to A$, $g(a)=ba$.
$g$ is right $R$-linear since
$$g(a\rightact r)=b(a\rightact r)=b\chi(ar)=\chi(bar)=(ba)\rightact r=g(a)r.$$
Finally
$$(p\circ g)(a)=p(ba)=ba+I=f(1_A)a=f(a).$$
\end{proof}

In \prref{2.1}, we characterize quasi-projectivity of $A$ as a right $R$-module.
Projectivity has been characterized in \cite[Prop. 2.4]{CVW}:

\begin{proposition}\prlabel{2.2}
Let $(R,i,\chi)$ be an $A$-ring with a right grouplike character. The following
assertions are equivalent.
\begin{enumerate}
\item $A$ is projective as a right $R$-module;
\item there exists $q\in Q$ such that $\chi(q)=1$.
\end{enumerate}
We refer to \cite[Prop. 2.4]{CVW} for more equivalent properties.
\end{proposition}

\begin{proposition}\prlabel{2.3}\cite[4.11]{Albu}
Let $R$ be a ring, $M$ a quasi-projective right $R$-module, and $N$ a
noetherian right $R$-module. Then $\Hom_R(M,N)$ is a noetherian right 
$\End_R(M)$-module.
\end{proposition}

\section{$I$-Frobenius rings}\selabel{3}
Let $(R,i)$ be an $A$-ring, and $I=(A,A,I,J,f,g)$ a strict Morita context
connecting $A$ with itself. We say that $R$ is an $I$-Frobenius $A$-ring if there exist
an element
$e=e^1\ot u^1\ot e^2\in R\ot_A I\ot_A R$ (summation
understood implicitely) and
an $A$-bimodule map
$\ol{\nu}:\ R\ot_A I\to A$ such that the following conditions are
satisfied, for all $r\in R$ and $u\in I$:
\begin{eqnarray}
&&re^1\ot u^1\ot e^2=e^1\ot u^1\ot e^2r;\eqlabel{3.1.1}\\
&&
\ol{\nu}(e^1\ot_A u^1)e^2=1_R;\eqlabel{3.1.2}\\
&&
e^1\ot_A u^1\ol{\nu}(e^2\ot_A u)=r1_R\ot_A u.\eqlabel{3.1.3}
\end{eqnarray}
If $I=(A,A,A,A,{\rm id}_A,{\rm id}_A)$, then the notion ``$I$-Frobenius" coincides with the
classical Frobenius property. Equivalent definitions are given in
\cite[Theorem 2.7]{CDM}.\\
$f:\ I\ot_A J\to A$ and $g:\ J\ot_A I\to A$ are $A$-bimodule isomorphisms, and
\begin{equation}\eqlabel{3.1.4}
f(u\ot_A v)u'=ug(v\ot_A u')~~;~~g(v\ot_A u)v'=vf(u\ot_A v'),
\end{equation}
for all $u,u'\in I$ and $v,v'\in J$. We will write
$$f^{-1}(1_A)=\sum_i u_i\ot v_i\in I\ot_A J.$$
From the fact that $f$ is an $A$-bimodule isomorphism, it follows easily that
\begin{equation}\eqlabel{3.1.5}
\sum_i au_i\ot v_i=\sum_i u_i\ot v_ia,
\end{equation}
for all $a\in A$. We have the following generalization of \cite[Theorem 2.7]{CVW}.

\begin{theorem}\thlabel{3.1}
Let $(R,i,\chi)$ be an $I$-Frobenius $A$-ring with a right grouplike character.
Then $J$ is an $(R,B)$-bimodule, with left
$R$-action
\begin{equation}\eqlabel{3.1.6}
r\cdot v=\sum_i \ol{\nu}(rg(v\ot \chi(e^1)u^1)e^2\ot_A u_i)v_i,
\end{equation}
and we have an isomorphism $\alpha:\ J\to Q$ of $(R,B)$-bimodules.
\end{theorem}

\begin{proof}
The map $\alpha$ is defined by the formula
\begin{equation}\eqlabel{3.1.7}
\alpha(v)=g(v\ot_A\chi(e^1)u^1)e^2,
\end{equation}
for all $v\in J$. Let us first show that $\alpha(v)\in Q$. For all $r\in R$,
we compute
\begin{eqnarray*}
&&\hspace*{-1cm}
\alpha(v)r= g(v\ot_A\chi(e^1)u^1)e^2r\equal{\equref{3.1.1}}
g(v\ot_A\chi(re^1)u^1)e^2\\
&\equal{\equref{1.1.0}}&g(v\ot_A\chi(\chi(r)e^1)u^1)e^2
\equal{\equref{3.1.1}}g(v\ot_A\chi(e^1)u^1)e^2\chi(r)=\alpha(v)\chi(r).
\end{eqnarray*}
$\alpha$ is right $B$-linear since
\begin{eqnarray*}
&&\hspace*{-2cm}\alpha(vb)=g(vb\ot_A\chi(e^1)u^1)e^2=
g(v\ot_Ab\chi(e^1)u^1)e^2\\
&=& g(v\ot_A\chi(be^1)u^1)e^2
g(v\ot_A\chi(e^1)u^1)e^2b\equal{\equref{3.1.1}}\alpha(v)b,
\end{eqnarray*}
for all $b\in B$. The inverse $\beta$ of $\alpha$ is given by the composition
$$Q\subset R\rTo^{R\ot_A f^{-1}}R\ot_AI\ot_AJ\rTo^{\ol{\nu}\ot_AJ}
A\ot_AJ\cong J,$$
or
$$\beta(q)=\sum_i\ol{\nu}(q\ot_A u_i)v_i,$$
for all $q\in Q$. Indeed, we compute for all $q\in Q$ that
\begin{eqnarray*}
&&\hspace*{-15mm}
\alpha(\beta(q))=
g(\sum_i\ol{\nu}(q\ot_A u_i)v_i\ot_A\chi(e^1)u^1)e^2\\
&=&\sum_ig(\ol{\nu}(q\ot_A u_i)v_i\chi(e^1)\ot_Au^1)e^2
\equal{\equref{3.1.5}}
\sum_ig(\ol{\nu}(q\ot_A \chi(e^1)u_i)v_i\ot_Au^1)e^2\\
&=&
\sum_ig(\ol{\nu}(q\chi(e^1)\ot_A u_i)v_i\ot_Au^1)e^2
\equal{\equref{3.1.1}}
\sum_ig(\ol{\nu}(\chi(e^1)\ot_A u_i)v_i\ot_Au^1)e^2q\\
&=&
\sum_i\ol{\nu}(\chi(e^1)\ot_A u_i) g(v_i\ot_Au^1)e^2q
= \sum_i\ol{\nu}(\chi(e^1)\ot_A u_i g(v_i\ot_Au^1))e^2q\\
&\equal{\equref{3.1.4}}&
\sum_i\ol{\nu}(\chi(e^1)\ot_A f(u_i \ot_Av_i)u^1)e^2q
= \ol{\nu}(e^1\ot_A e^1)e^2q=q.
\end{eqnarray*}
For all $v\in J$, we have that
\begin{eqnarray*}
&&\hspace*{-15mm}
\beta(\alpha(v))=
\sum_i\ol{\nu}(g(v\ot_A\chi(e^1)u^1)e^2\ot_A u_i)v_i=
\sum_ig(v\ot_A\chi(e^1)u^1)\ol{\nu}(e^2\ot_A u_i)v_i\\
&=&\sum_ig(v\ot_A\chi(e^1)u^1\ol{\nu}(e^2\ot_A u_i))v_i
\equal{\equref{3.1.3}} \sum_i g(v\ot_A \chi(1_R)u_i)v_i\\
&=&\sum_i g(v\ot_A u_i)v_i
\equal{\equref{3.1.4}} \sum_i vf(u_i\ot_A v_i)=v.
\end{eqnarray*}
This shows that $\alpha$ is an isomorphism of right $B$-modules. We can transport
the left $B$-action on $Q$ to $J$ such that $\alpha$ becomes an $(R,B)$-bimodule
map. This yields formula \equref{3.1.6}.
\end{proof}

The composition
$$\tr=\chi\circ \alpha:\ J\to Q\to B$$
is a $B$-bimodule map (see \leref{1.1}), and will be called the {\sl trace map}.
It is given by the formula
\begin{equation}\eqlabel{3.18}
\tr(v)=\chi(g(v\ot_A\chi(e^1)u^1)e^2).
\end{equation}
Combining \prref{2.2} and \thref{3.1}, we obtain the following result:

\begin{proposition}\prlabel{3.2}
Let $(R,i,\chi)$ be an $I$-Frobenius $A$-ring with a right grouplike character.
The following assertions are equivalent.
\begin{enumerate}
\item $A$ is projective as a right $R$-module;
\item there exists $v\in J$ such that $\tr(v)=1_B$.
\end{enumerate}
\end{proposition}

Now assume that $R$ is Frobenius $A$-ring, that is, $I=A$. Then the above formulas
simplify. $e=e^1\ot e^2\in R\ot_AR$, $\ol{\nu}:\ R\to A$ is an $A$-bimodule
map, and the trace map $\tr:\ A\to B$ is given by
$$\tr(a)= \chi(a\chi(e^1)e^2).$$

\section{Fully bounded noetherian rings}\selabel{4}
We recall some definitions and basic results from \cite{Garcia}.
Let $R$ be a ring, and $M,P\in \Mm_R$. For a subset $X$ of $\Hom_R(P,M)$, we
write
$$r_P(X)=\cap\{\Ker f~|~f\in X\}.$$
In particular, for $X\subset M\cong \Hom_R(R,M)$, we have
$$r_R(X)=\{r\in R~|~xr=0\}.$$
$M$ is called {\sl finitely} $P$-{\sl generated} if there exists an epimorphism of right
$R$-modules $P^n\to M\to 0$.\\
$M$ is called $P$-{\sl faithful} if $\Hom_R(P,M')\neq 0$, for every nonzero 
submodule $M'\subset M$.\\
$R$ is called {\sl right bounded} if every essential right ideal contains a non-zero
two-sided ideal. $R$ is called {\sl right fully bounded} if $R/P$ is right bounded,
for every two-sided prime ideal $P$ of $R$. A ring $R$ that is right fully bounded
and right noetherian is called a {\sl right fully bounded noetherian ring} or a
{\sl FBN ring}.  Characterizations of right FBN rings are given in \cite[Theorem 1.2]{Garcia}.
For later use, we recall one of them.

\begin{proposition}\prlabel{4.1}
For a ring $R$, the following conditions are equivalent.
\begin{enumerate}
\item $R$ is right FBN;
\item for every finitely generated right $R$-module $M$, there exists a finite subset $F\subset M$
such that $r_R(M)=r_R(F)$.
\end{enumerate}
\end{proposition}

A right $R$-module $P$ is called a {\sl right FBN-module} if it is noetherian and for
every finitely generated $P$-faithful right $R$-module $M$, there exists a finite
subset $F\subset \Hom_R(P,M)$ such that $r_P(F)=r_P(\Hom_R(P,M))$. We recall
the following properties from \cite{Garcia}.

\begin{proposition}\prlabel{4.2} \cite[Theorem 1.7]{Garcia}
For a quasi-projective, noetherian right $R$-module $P$, the following assertions are
equivalent:
\begin{enumerate}
\item $\End_R(P)$ is right FBN;
\item $P$ is an FBN right $R$-module.
\end{enumerate}
\end{proposition}

\begin{proposition}\prlabel{4.3} \cite[Corollary 1.8]{Garcia}
Let $P$ be a quasi-projective FBN right $R$-module, $Q$ a finitely $P$-generated
right $R$-module, and $M$ a finitely generated $Q$-faithful right $R$-module. 
For every $X\subset \Hom_R(Q,M)$, there exists a finite subset
$F\subset X$ such that $r_Q(X)=r_Q(F)$.
\end{proposition}

\begin{proposition}\prlabel{4.4} \cite[Corollary 1.9]{Garcia}
A right noetherian ring $R$ is right FBN if and only if every finitely generated right
$R$-module is FBN.
\end{proposition}

We can now state the main result of this paper.

\begin{theorem}\thlabel{4.5}
Let $(R,i,\chi)$ be an $A$-ring with a right grouplike character, and consider the
following statements.
\begin{enumerate}
\item $R\in \Mm_A$ is finitely generated and $A$ is right FBN;
\item $R$ is right FBN and $A$ is right noetherian;
\item $B$ is right FBN and $A$ is right noetherian.
\end{enumerate}
Then $1)\Rightarrow 2)$.\\
If $A$ is quasi-projective as a right $R$-module, then $2)\Rightarrow 3)$.\\
If $A$ is projective as a right $R$-module and $R$ is an $I$-Frobenius $A$-ring for
some strict Morita context $I=(A,A,I,J,f,g)$, then $3)\Rightarrow 1)$ and the three
conditions are equivalent.
\end{theorem}

\begin{proof}
$1)\Rightarrow 2)$. It follows from \prref{4.4} that $R$ is an FBN right $R$-module.
Let $M$ be a finitely generated right $R$-module; then $M$ is also finitely generated
as a right $A$-module. We claim that $M$ is an $R$-faithful right $A$-module.
Indeed, take a non-zero right $A$-module $M'\subset M$. Since $M'\cong
\Hom_A(A,M')$, there exists a non-zero $f\in \Hom_A(A,M')$, and the composition
$f\circ \chi:\ R\to M'$ is non-zero, since $\chi$ is surjective.\\
Now take $P=R$, $Q=A$ in \prref{4.3}, and consider the subset
$M\cong \Hom_R(R,M)\subset \Hom_A(R,M)$. It follows that there exists a finite
$F\subset M$ such that $r_A(F)=r_A(M)$. It then follows from \prref{4.1} that $R$
is right FBN.\\

$2)\Rightarrow 3)$. $A$ is a finitely generated (even cyclic) right $R$-module, so
it follows from \prref{4.4} that $A$ is an FBN right $R$-module. It then follows from
\prref{4.2} that $\End_R(A)\cong B$ is right FBN.\\

$3)\Rightarrow 1)$. We will apply \prref{2.3} with $M=A$ and $N=R$.
By assumption, $A$ is quasi-projective as a right $R$-module. Since $R/A$ is
$I$-Frobenius, $R$ is finitely generated projective as a right $R$-module. Since
$A$ is right noetherian, $R$ is also right noetherian.\\
It follows from \leref{1.2}, \prref{2.3} and \thref{3.1} that $\Hom_R(A,R)\cong R^R=Q\cong J$
is noetherian as a right module over $\End_R(A)\cong A^R=B$. It then follows that $J$
is finitely generated as a right $B$-module. Let $\{e_1,\cdots,e_k\}$ be a set of
generators of $J$ as a right $B$-module.\\
Recall that we have an $A$-bimodule isomorphism $f:\ I\ot_A J\to A$. With notation
as in \seref{3}, we have, for $a\in A$,
$$f^{-1}(a)=\sum_{i=1}^n u_i\ot_A v_ia\in I\ot_A J.$$
For every $i$, we can find $b_{i1},\cdots,b_{ik_i}\in B$ such that
$$v_ia=\sum_{j=1}^{k_i}e_jb_{ij}.$$
We then easily compute that
\begin{eqnarray*}
a&=& f\Bigl(\sum_{i=1}^n u_i\ot_A v_ia\Bigr)=
f\Bigl(\sum_{i=1}^n\sum_{j=1}^{k_i} u_i\ot_A e_jb_{ij}\Bigr)=\
\sum_{i=1}^n\sum_{j=1}^{k_i} f(u_i\ot_Ae_j)b_{ij},
\end{eqnarray*}
and we conclude that $A$ is finitely generated as a right $B$-module.\\
Take $M\in \Mm_A$ finitely generated. Then $M$ is also finitely generated as a right
$B$-module. We now show that $M$ is an $A$-faithful right $B$-module.
Let $M'$ be a non-zero right $B$-submodule of $M$, and take $0\neq m'\in M'$.
It follows from \prref{3.2} that there exists $v\in J$ such that $\tr(v)=1_B$. The map
$f:\ A\to M$, $f(a)=m'\tr(va)$ is right $B$-linear, and different from $0$ since
$f(1_A)=m'\neq 0$.\\
Observe now that
\begin{itemize}
\item $B$ is a quasi-projective FBN right $B$-module;
\item $A$ is a finitely $B$-generated right $B$-module;
\item $M$ is a finitely generated $A$-faithful right $B$-module.
\end{itemize}
Applying \prref{4.3} to $M\cong \Hom_A(A,M)\subset \Hom_B(A,M)$, we find that
there exists a finite subset $F\subset M$ such that $r_A(F)= r_A(M)$.
It then follows from \prref{4.1} that $A$ is right FBN.
\end{proof}

\begin{remark}\relabel{4.6}
We do not know whether the implication $3)\Rightarrow 1)$ holds under the weaker
assumption that $A\in \Mm_R$ is quasi-projective. The projectivity is used at the point
where we applied \prref{4.3}.
\end{remark}

\section{Application to Frobenius algebras}\selabel{5}
Let $k$ be a commutative ring, and consider two $k$-algebras $A$ and $S$,
and an algebra map $j:\ S\to A$. All unadorned tensor products in this Section are over $k$. It is easy to establish that
$(R=S^{\rm op}\ot A, i,\chi)$ with
$$i:\ A\to S^{\rm op} \ot A,~~i(a)=1_S\ot a,$$
$$\chi:\ S^{\rm op}\ot A\to A,~~\chi(s\ot a)=j(s)a$$
is an $A$-ring with a right grouplike character. Also observe that
the categories $\Mm_R$ and ${}_S\Mm_A$ are isomorphic. For $M\in {}_S\Mm_A$,
we have that
$$M^R=\{m\in M~|~sm=mj(s),~{\rm for~all~}s\in S\}=C_S(M).$$
In particular, $B=A^R=C_S(A)$ and
$$Q=\{\sum_i s_i\ot a_i\in S^{\rm op} \ot A~|~
\sum_i ts_i\ot a_i=\sum_i s_i\ot a_ij(t),~{\rm for~all~}t\in S\}.$$
Consequently $A$ is projective as a right $R$-module if and only if there exists
$\sum_i s_i\ot a_i\in Q$ such that $\sum_i j(s_i)a_i=1_A$.\\
From \prref{2.1}, it follows that $A$ is quasi-projective as a right $R$-module
if and only if for every $(S,A)$-submodule $I$ of $A$ and $a\in A$ such that
$as-sa\in I$, for all $s\in S$, there exists $b\in B$ such that $a-b\in I$.\\
Assume that $S$ is a Frobenius $k$-algebra, with Frobenius system
$(e=e^1\ot e^2,\ol{\nu})$. Then $S^{\rm op}$ is also a Frobenius algebra, with
Frobenius system $(e=e^2\ot e^1,\ol{\nu})$, and $S^{\rm op}\ot A$ is a Frobenius
$A$-ring, with Frobenius system $(E, N)$, with
$E=(e^2\ot 1_A)\ot_A (e^1\ot 1_A)$ and
$$N:\ S^{\rm op}\ot A\to A,~~N(s\ot a)=\ol{\nu}(s)a.$$
We then have the isomorphism
$$\alpha:\ A\to Q,~~\alpha(a)=e^1\ot aj(e^2)$$
and the trace map
$$\tr:\ A\to B,~~\tr(a)=j(e^1)aj(e^2).$$
$A$ is projective as a right $R$-module if and only if there exists $a\in A$ such that
$\tr(a)=1$.

\begin{corollary}\colabel{5.1}
Let $S$ be a Frobenius algebra over a commutative ring $k$, and $j:\ S\to A$
an algebra map. Furthermore, assume that 
there exists $a\in A$ such that $\tr(a)=1$. Then the following assertions are
equivalent:
\begin{enumerate}
\item $A$ is right FBN;
\item $S^{\rm op}\ot A$ is right FBN and $A$ is right noetherian;
\item $B=C_S(A)$ is right FBN and $A$ is right noetherian.
\end{enumerate}
\end{corollary}

\section{Application to Hopf algeba actions}\selabel{6}
Let $H$ be a finitely generated projective Hopf algebra over a commutative ring
$k$, and $A$ a left $H$-module algebra. The smash product $R=A\# H$ is equal to $A\ot H$ as
a $k$-module, with multiplication given by the formula
$$(a\# h)(b\# k)=a(h_{(1)}\cdot b)\# h_{(2)}k.$$
The unit is $1_A\# 1_H$. Consider the maps
$$i:\ A\to A\#H,~~i(a)=a\#1_H,$$
$$\chi:\ A\# H\to A,~~\chi(a\# h)=a\varepsilon(h).$$
Straightforward computations show that $(A\# H,i,\chi)$ is an $A$-ring
with a left grouplike character. It is also easy to prove that
$$A^R=\{a\in A~|~h\cdot a=\varepsilon(h)a,~{\rm for~all~}h\in H\}=A^H$$
is the subalgebra of invariants of $R$.\\
In a similar way, we can associate an $A$-ring with right grouplike character to
a right $H$-comodule algebra. We will discuss the left handed case here, in order
to recover the results from \cite{Sorin,Garcia,Guedenon}. The results from the previous
Sections can easily be restated for rings with a left grouplike character.\\
Let $I=\int_{H^*}^l$ and $J=\int_H^l$ be the spaces of left integrals on and
in $H$. $I$ and $J$ are projective rank one $k$-modules, and $H/k$ is
$I$-Frobenius (see for example \cite[Theorem 3.4]{CDM}). We need an explicit
description of the Frobenius system. From the Fundamental Theorem, it follows
that we have an isomorphism
$$\phi:\ I\ot H\to H^*,~~\phi(\varphi\ot h)=h\cdot \varphi,$$
with $\lan h\cdot \varphi,k\ran=\lan \varphi,kh\ran$. If $t\in J$, then
$$\phi(\varphi\ot t)(h)=\lan \varphi,ht\ran=\lan\varphi,t\ran\varepsilon(h),$$
so $\phi$ restricts to a monomorphism $\tilde{\phi}:\ I\ot J\to k\varepsilon$.
If $I$ and $J$ are free of rank one, then $\tilde{\phi}$ is an isomorphism,
as there exist $\varphi\in I$ and $t\in J$ such that $\lan\varphi,t\ran=1$
(see for example \cite[Theorem 31]{CMZ}, \cite{Pareigis71}. Hence $\tilde{\phi}$ is an isomorphism
after we localize at a prime ideal $p$ of $k$, and this implies that $\tilde{\phi}$
is itself an isomorphism. Consequently $J^*\cong I$. Consider
$\tilde{\phi}^{-1}(\varepsilon)=\sum_i\varphi_i\ot t_i\in I\ot J$. Then
\begin{equation}\eqlabel{6.1.1}
\sum_i \lan\varphi_i,t_i\ran=1.
\end{equation}
Furthermore $\{(\varphi_i,t_i)~|~i=1,\cdots,n\}$ is a finite dual basis for $I$,
so we have $t=\sum_i \lan\varphi_i,t\ran t_i$, $\varphi=\sum_i\lan\varphi,t_i\ran
\varphi$, for all $t\in J$ and $\varphi\in I$. $\phi$ induces an isomorphism
$$\psi:\ H\to H^*\ot J,~~\psi(h)=\sum_i h\cdot\varphi_i\ot t_i.$$
The inverse of $\psi$ is given by the formula
$$\psi^{-1}(h^*\ot t)=\lan h^*,\ol{S}(t_{(1)})\ran t_{(2)},$$
where $\ol{S}$ is the inverse of the antipode $S$; recall from
\cite{Pareigis71} that the antipode of
a finitely generated projective Hopf algebra is always bijective. Indeed,
it is straightforward to show that $\psi^{-1}$ is a right inverse of $\psi$.
First observe that
$$\psi(\psi^{-1}(h^*\ot t))=\sum_i\lan h^*,\ol{S}(t_{(1)})\ran t_{(2)}\cdot \varphi_i\ot t_i.$$
Now we compute for all $h\in H$ that
\begin{eqnarray*}
&&\hspace*{-2cm}
\lan h^*,\ol{S}(t_{(1)})\ran \lan t_{(2)}\cdot\varphi_i,h\ran=
\lan h^*,\ol{S}(t_{(1)})\ol{S}(h_{(2)})h_{(1)}\ran \lan  \varphi_i,h_{(3)}t_{(2)}\ran\\
&=& \lan h^*,\ol{S}(h_{(2)}t_{(1)})h_{(1)}\ran \lan  \varphi_i,h_{(3)}t_{(2)}\ran\\
&=& \lan h^*,\ol{S}(1_H)h_{(1)}\ran \lan  \varphi_i,h_{(2)}t\ran
= \lan h^*,h\ran \lan  \varphi_i,t\ran,
\end{eqnarray*}
where we used the fact that $\varphi_i$ and $t$ are integrals. It follows that
$$\psi(\psi^{-1}(h^*\ot t))=\sum_i h^*\ot \lan\varphi_i,t\ran t_i=h^*\ot t.$$
A right inverse of an invertible element is also a left inverse, so it follows that
$$
1_H=\psi(\psi^{-1}(1_H))=\sum_i\lan \varphi_i,\ol{S}(t_{i(1)})\ran t_{i(2)}=
 \sum_i\lan \varphi_i\circ \ol{S},t_{i(1)}\ran t_{i(2)}=
\sum_i\lan \varphi_i\circ \ol{S},t_i\ran 1_H,$$
where we used the fact that $\varphi_i\circ \ol{S}$ is a right integral on $H$.
We conclude that
\begin{equation}\eqlabel{6.1.2}
\sum_i\lan \varphi_i,\ol{S}(t_i)\ran=1.
\end{equation}
Consider the particular situation where $I$ and $J$ are free rank one modules.
Then there exist free generators
$\varphi_1$ of $ I$ and $t_1$ of $ J$ such that $\lan\varphi_1,t_1\ran=1$.
From \equref{6.1.2} it follows that $\lan\varphi_1,\ol{S}(t_1)\ran=1$. 
For arbitrary $\varphi=x\varphi_1\in I$ and $t=yt_1\in J$, it then follows that
$\lan \varphi,t\ran=xy \lan\varphi_1,t_1\ran=xy=
xy\lan\varphi_1,\ol{S}(t_1)\ran=\lan\varphi,\ol{S}(t)\ran$.
Consider the case where $I$ and $J$ are not necessarily free, and take
$\varphi\in I$, $t\in J$ and a prime ideal $p$ of $k$. 
Then the images of $\lan \varphi,t\ran$ and
$\lan\varphi,\ol{S}(t)\ran$ in the localized ring $k_p$ are equal, since the integral
space of the Hopf $k_p$-algebra $H_p$ is free. So we can conclude that
\begin{equation}\eqlabel{6.1.3}
\lan \varphi,t\ran=\lan\varphi,\ol{S}(t)\ran.
\end{equation}

\begin{lemma}\lelabel{6.1}
Let $H$ be a finitely generated projective Hopf algebra over a commutative ring
$k$. There exist $t_i\in J=\int_H^l$ and $\varphi_i\in I=\int_{H^*}^l$ such that
$\sum_i \lan\varphi_i,t_i\ran=1$. $H$ is an $I$-Frobenius $k$-algebra, with Frobenius system
$(e,\ol{\nu})$ with
\begin{eqnarray*}
&&e=\sum_i t_{i(2)}\ot\varphi_i\ot\ol{S}(t_{i(1)})\\
&&\ol{\nu}=\sum_j t_j\ot \varphi_j\in (H\ot I)^*\cong J\ot H^*
\end{eqnarray*}
\end{lemma}

\begin{proof}
It is straightforward to show that $e\in C_H(H\ot I\ot H)$; this also follows
from \cite[Prop. 3.3]{CDM}, taking into account that $e=
i'(\varphi\ot \ol{S}(t))$.\\
Write $e=e^1\ot u^1\ot e^2\in H\ot I\ot H$. We compute that
\begin{eqnarray}
&&\hspace*{-2cm}
\ol{\nu}(e^1\ot u^1\ot e^2)=
\sum_{i,j}\lan \varphi_j,t_{i(2)}\ran \lan \varphi,t_j\ran \ol{S}(t_{i(1)})\nonumber\\
&=& \sum_{i}\lan \varphi_i,t_{i(2)}\ran \ol{S}(t_{i(1)})
= \sum_i\ol{S}(\lan \varphi_i,t_i\ran 1_H)\equal{\equref{6.1.1}}1_H.\eqlabel{6.1.4}
\end{eqnarray}
For all $\varphi\in I$, we calculate
\begin{eqnarray}
&&\hspace*{-2cm}
e^1\ot u^1\ol{\nu}( e^2\ot \varphi)=
\sum_{i,j} t_{i(2)}\ot\varphi_i\lan\varphi_j,\ol{S}(t_{i(1)})\ran\lan\varphi,t_j\ran\nonumber\\
&=& \sum_{i,j} 1_H\ot\varphi_i\lan\varphi_j,\ol{S}(t_{i})\ran\lan\varphi,t_j\ran
= \sum_{i} 1_H\ot\varphi_i\lan\varphi,\ol{S}(t_{i})\ran\nonumber\\
&\equal{\equref{6.1.3}}& \sum_{i} 1_H\ot\varphi_i\lan\varphi,t_{i}\ran=1_H\ot \varphi.
\eqlabel{6.1.5}
\end{eqnarray}
It now follows from \cite[Theorem 3.1]{CDM} that $(e,\ol{\nu})$ is a
Frobenius system.
\end{proof}

\begin{proposition}\prlabel{6.2}
Let $H$ be a finitely generated projective Hopf algebra over a commutative ring
$k$, and $A$ a left $H$-module algebra. Then $A\ot H$ is an $A\ot I$-Frobenius
$A$-algebra,
with Frobenius system $(E,N)$, with
\begin{eqnarray*}
&&\hspace*{-2cm}
E=E^1\ot_AU^1\ot_AE^2=(1_A\# e^1)\ot_A(1_A\ot u^1)\ot_A(1_A\# e^1)\\
&=&\sum_i (1_A\# t_{i(2)})\ot_A(1_A\ot\varphi_i)
\ot_A(1_A\#\ol{S}(t_{i(1)})),\\
&&\hspace*{-2cm}N:\ (A\#H)\ot_A (A\ot I)\cong A\# H\ot I\to A,\\
&&N(a\#h\ot\varphi)=a\ol{\nu}(h\ot\varphi)=\sum_j a\lan\varphi_j,h\ran\lan \varphi,t_j\ran.
\end{eqnarray*}
Here we used the notation introduced above.
\end{proposition}

\begin{proof}
The proof is an adaptation of the proof of \cite[Proposition 5.1]{CVW}. Let us first show
that $E$ satisfies \equref{3.1.1}.
\begin{eqnarray*}
&&\hspace*{-15mm}
\sum_i (1_A\# t_{i(2)})\ot_A(1_A\ot\varphi_i)\ot_A(1_A\#\ol{S}(t_{i(1)})(a\# h)\\
&=&
\sum_i (1_A\# t_{i(3)})\ot_A(1_A\ot\varphi_i)\ot_A(\ol{S}(t_{i(2)})\cdot a \# \ol{S}(t_{i(1)}) h)\\
&=&
\sum_i (1_A\# t_{i(3)})\ot_A(\ol{S}(t_{i(2)})\cdot a\ot\varphi_i)\ot_A(1_A \# \ol{S}(t_{i(1)}) h)\\
&=&
\sum_i ((t_{i(3)}\ol{S}(t_{i(2)}))\cdot a\# t_{i(4)})\ot_A(1_A\ot\varphi)\ot_A(1_A \# \ol{S}(t_{i(1)}) h)\\
&=&
\sum_i ( a\# t_{i(2)})\ot_A(1_A\ot\varphi)\ot_A(1_A \# \ol{S}(t_{i(1)}) h)\\
&=&
\sum_i ( a\#h t_{i(2)})\ot_A(1_A\ot\varphi)\ot_A(1_A \# \ol{S}(t_{i(1)}) )\\
&=&
\sum_i ( a\#h)(1_A\# t_{i(2)})\ot_A(1_A\ot\varphi_i)\ot_A(1_A\#\ol{S}(t_{i(1)}).
\end{eqnarray*}
Obviously $N$ is left $A$-linear. Right $A$-linearity can be proved as follows:
\begin{eqnarray*}
&&\hspace*{-2cm}
N((1\# h\ot \varphi)a)=N(h_{(1)}a\# h_{(2)}\ot\varphi)\\
&=& \sum_j h_{(1)}\cdot a\lan\varphi_j,h_{(2)}\ran\lan \varphi,t_j\ran
= N(1\# h\ot \varphi)a.
\end{eqnarray*}
\equref{3.1.2} is satisfied since
\begin{eqnarray*}
&&\hspace*{-2cm}
N(E^1\ot_AU^1)E^2=1_A\ol{\nu}(e^1\ot u^1)(1_A\# e^2)\\
&\equal{\equref{6.1.4}}& 1_A\#\ol{\nu}(e^1\ot u^1) e^2=1_A\#1_H.
\end{eqnarray*}
Let us finally show that \equref{3.1.3} holds. For all $a\in A$ and $\varphi\in I$, we have
\begin{eqnarray*}
&&\hspace*{-2cm}
E^1\ot_AU^1N(E^2\ot_A(a\ot\varphi))\\
&=&
\sum_i(1_A\#t_{i(2)})\ot_A(1_A\ot\varphi_i)N(a\#\ol{S}(t_{i(1)})\ot\varphi)\\
&=&\sum_{i,j}(1_A\#t_{i(2)})\ot_A(a\ot\varphi_i)\lan\varphi_j,\ol{S}(t_{i(1)})
\lan\varphi,t_j\ran\\
&\equal{\equref{6.1.4}}& (1_A\#1_H)\ot_A (a\ot\varphi)
\end{eqnarray*}
\end{proof}

\begin{proposition}\prlabel{6.3}
Let $H$ be a finitely generated projective Hopf algebra, and $A$ a left $H$-module
algebra. The trace map $\tr:\ A\ot J\to B=A^H$ is given by the formula
$$\tr(a\ot t)=t\cdot a.$$
\end{proposition}

\begin{proof}
Observe that the map $g:\ (J\ot A)\ot_A(I\ot A)$ in the Morita context associated to
$I\ot A$ is given by the formula
$$g((t\ot a)\ot_A (\varphi\ot b))=\lan \varphi,t\ran ab.$$
Using the left handed version of \equref{3.18}, we compute, for $V=a\ot t\in A\ot J$
that
\begin{eqnarray*}
&&\hspace*{-15mm}
\tr(a\ot v)=\chi(E^1g(U^1\chi(E^2)\ot_AV))
=\sum_i \chi((1_A\# t_i)g((1_A\ot\varphi)\ot (a\ot t)))\\
&=&\sum_i \chi((1_A\# t_i)a\lan\varphi,t\ran)=\chi((1_A\# t)a)
=\chi(t_{(1)}\cdot a\# t_{(2)})=t\cdot a.
\end{eqnarray*}
\end{proof}

We can now apply Propositions \ref{pr:2.1}, \ref{pr:2.2} and \ref{pr:3.2}, and \thref{4.5},
and obtain the following result.

\begin{corollary}\colabel{6.4}
Let $H$ be a finitely generated projective Hopf algebra, and $A$ a left $H$-module
algebra. Assume that there exist $a_i\in A$ and $t_i\in \int_l^H$
such that $\sum_it_i\cdot a_i=1$.\\
Then the following assertions are equivalent;
\begin{enumerate}
\item $A$ is left FBN;
\item $A\# H$ is left FBN and $A$ is left noetherian;
\item $B$ is left FBN and $A$ is left noetherian.
\end{enumerate}
\end{corollary}

We recover \cite[Theorem 2.3 and Corollary 2.4]{Garcia}, \cite[Theorem 8]{Sorin}
and \cite[Theorem 2.4]{Guedenon}. If $H$ is Frobenius (e.g. if $k$ is a field, or
$H=kG$ is a finite group algebra), then the space of left integrals is free. We can then
take a free generator $t$ of $\int_H^l$ and the condition of the trace map means that
there exists $a\in A$ such that $t\cdot a=1$. We observe that - in the case where
the space of integrals is not free - the sufficient condition in \coref{6.4} that there exist $a_i\in A$ and $t_i\in \int_l^H$
such that $\sum_it_i\cdot a_i=1$ is weaker than the one given in \cite[Theorem 8]{Sorin},
where a single $t\in \int_l^H$ and $a\in A$ with $t\cdot a=1$ are needed.\\
In \cite{Garcia} and \cite{Guedenon}, it is stated that \coref{6.4} holds under the weaker
assumption (called (C1)) that $A$ is $A\# H$-{\sl quasi}-projective. There seems to be
a hole in the proofs in \cite{Garcia} and \cite{Guedenon}: the proof of the implication $3)\Longrightarrow 1)$ uses the projectivity of $A$ as an $A\# H$-module (see
\reref{4.6}).

\section{Application to corings}\selabel{7}
Let $A$ be a ring. An $A$-coring is a coalgebra in the category of
$A$-bimodules ${}_A\Mm_A$. This means that we have two $A$-bimodule maps
$$\Delta_\Cc:\ \Cc\to \Cc\ot_A\Cc~~{\rm and}~~\varepsilon_\Cc:\ \Cc\to A$$
satisfying some coassociativity and counit axioms. The maps $\Delta_\Cc$
and $\varepsilon_\Cc$ are called the comultiplication and counit, and we
use the Sweedler notation
$$\Delta_\Cc(c)=c_{(1)}\ot_A c_{(2)},$$
where summation is understood implicitely. Corings were revived recently
in \cite{Brzezinski02}, and we refer to \cite{BrzezinskiWisbauer} for a
detailed discussion of all kinds of applications. The left dual
$R={}^*\Cc={}_A\Hom(\Cc,A)$ is an $A$-ring, with multiplication rule
$$(f\# g)(c)=g(c_{(1)}f(c_{(2)})),$$
for all $c\in \Cc$ and $f,g\in {}^*\Cc$. The unit is $\varepsilon_\Cc$,
and the ring morphism $i:\ A\to{}^* \Cc$ is given by
$$i(a)(c)=\varepsilon_\Cc(c)a.$$
The $A$-bimodule structure on ${}^*\Cc$ is then given by the formula
$$(afb)(c)=f(ca)b,$$
for all $a,b\in A$, $f\in \Cc^*$ and $c\in \Cc$.\\
$x\in \Cc$ is called grouplike if $\Delta_\Cc(x)=x\ot_R x$ and
$\varepsilon_\Cc(x)=1_R$. $(\Cc,x)$ is then called an $R$-coring with
a grouplike element. Now consider the map
$\chi:\ {}^*\Cc\to A$, $\chi(f)=f(x)$.
It can be shown easily (see \cite{CVW}) that $(\Cc^*,i,\chi)$ is
an $A$-ring with a right grouplike character.
We can also compute that
$$B=A^R=\{a\in A~|~f(xa)=af(x),~{\rm for~all~}f\in \Cc^*\}.$$
Using the grouplike element $x$, we can define a right $\Cc$-coaction on $A$,
namely
$$\rho:\ \Cc\to A\ot_A\Cc\cong \Cc,~~\rho(r)=1_A\ot_A xa=xa.$$
We can consider the subring of coinvariants
$$A^{{\rm co}\Cc}=\{a\in A~|~ax=xa\}.$$
In general, $A^{{\rm co}\Cc}$ is a subring of $A^R$, and they are equal
if $\Cc$ is finitely generated and projective as a right $A$-module.\\
An $A$-coring $\Cc$ is called Frobenius if there exist an $A$-bimodule
map $\theta:\ \Cc\ot_A\Cc\to A$ and $z\in C_A( \Cc)$ (that is,
$az=za$, for all $a\in A$) such that the following
conditions hold, for all $c,d\in \Cc$:
$$c_{(1)}\theta(c_{(2)}\ot d)=\theta(c\ot d_{(1)})d_{(2)},$$
$$\theta(z\ot c)=\theta(c\ot z)=1.$$
We refer to \cite[Theorem 35]{CMZ} for the explanation of this definition.
If $\Cc$ is Frobenius, $\Cc$ is finitely generated and projective as a
(left or right) $A$-module, and ${}^*\Cc/A$ is Frobenius (see
\cite[Theorem 36]{CMZ}). Then we also have (see \cite[Sec. 3]{CVW}) that
$$Q=\{q\in {}^*\Cc~|~c_{(1)}q(c_{(2)})=q(c)x\}.$$
It follows from \cite[Theorem 2.7]{CVW} or \thref{3.1} that we have an isomorphism
of $({}^*\Cc,B)$-bimodules $\alpha:\ A\to Q$, given by
$$\alpha(a)(c)=\theta(ca\ot_A x),$$
for all $a\in A$ and $c\in \Cc$.
The inverse $\alpha^{-1}$ is given by $\alpha^{-1}(q)=q(z)$, and the left
${}^*\Cc$-action on $A$ is 
$$f\cdot a=\theta(z_{(1)}f(z_{(2)})\ot_A x).$$
This can be verified directly as follows:
$$\alpha(\alpha^{-1}(a))=\theta(za\ot_A x)=\theta(az\ot x)=a\theta(z\ot x)=a,$$
and
\begin{eqnarray*}
&&\hspace*{-2cm}
\alpha(\alpha^{-1}(q))(c)=\theta(cq(z)\ot_A x)=\theta(c\ot_A q(z)x)
= \theta(c\ot_A z_{(1)}q(z_{(2)})\\
&=& \theta(c\ot_A z_{(1)})q(z_{(2)})
=q(\theta(c\ot_A z_{(1)})z_{(2)})\\
&=&q(c_{(1)}\theta(c_{(2)}\ot z))
= q(c_{(1)}\varepsilon(c_{(2)}))=q(c).
\end{eqnarray*}
The trace map $\tr:\ A\to B$ is given by
$$\tr(a)=\theta(xa\ot_A x).$$

\begin{corollary}\colabel{7.1}
Let $(\Cc,x)$ be a Frobenius $A$-coring with a fixed grouplike element, and
Frobenius system $(\theta, z)$, and assume that there exists $a\in A$ such that
$\tr(a)=1$. Then the following assertions are equivalent.
\begin{enumerate}
\item $A$ is right FBN;
\item ${}^*\Cc$ is right FBN and $A$ is right noetherian;
\item $B=A^{{\rm co}\Cc}$ is right FBN and $A$ is right noetherian.
\end{enumerate}
\end{corollary}

\begin{center}
{\bf Acknowledgement}
\end{center}
We thank Angel del R\'{\i}o and Sorin D\v{a}sc\v{a}lescu for discussing with us
the sufficiency of the quasi-projectivity assumption in the proof of $3)\Longrightarrow 1)$
in \thref{4.5}.


\begin{thebibliography}{99}

\bibitem{Albu}
T. Albu and C. N\v{a}st\v{a}sescu, ``Relative finiteness in module 
theory", {\sl Monographs Textbooks Pure Appl. Math.} {\bf 84}, Marcel Dekker, New York, 1984.

\bibitem{Brzezinski02}
T. Brzezi\'nski, {The structure of corings. Induction functors,
Maschke-type theorem, and Frobenius and Galois properties},
{\sl Algebr. Representat. Theory} \textbf{5} (2002), 389--410.

\bibitem{BrzezinskiWisbauer}
T. Brzezi\'nski and R. Wisbauer, ``Corings and comodules", 
{\sl London Math. Soc. Lect. Note Ser.} {\bf 309},
  Cambridge University Press, Cambridge, 2003.

\bibitem{CDM}
S. Caenepeel, E. De Groot and G. Militaru, Frobenius functors of the
second kind, {\sl Comm. Algebra} {\bf 30} (2002), 5359--5391.

\bibitem{CMZ}
S. Caenepeel, G. Militaru and Zhu Shenglin ``Frobenius and separable 
functors for generalized module categories and nonlinear equations", 
{\sl Lect. Notes in Math.} \textbf{1787}, Springer Verlag, Berlin, 2002.

\bibitem{CVW}
S. Caenepeel, J. Vercruysse and Shuanhong Wang,
{Morita Theory for corings and cleft entwining structures},
{\sl J. Algebra} {\bf 276} (2004), 210--235.

\bibitem{Sorin}
S. D\v{a}sc\v{a}lescu, A. Kelarev, and B. Torrecillas: FBN
Hopf module algebras, {\sl Comm. Algebra} {\bf 25} (1997), 3521--3529.

\bibitem{Kaoutit}
L. El Kaoutit, J. G\'omez-Torrecillas and F. Lobillo,
Semisimple corings, {\sl Algebra Coll.} {\bf 11} (2004), 427--442.

\bibitem{Garcia} 
J. J. Garc\'{\i}a and A. Del R\'{\i}o, Actions of groups on fully bounded
Noetherian rings, {\sl Comm. Algebra} {\bf 22} (1994), 1495--1505. 

\bibitem{Guedenon} 
T. Gu\'ed\'enon, Actions of Hopf algebras on fully bounded
Noetherian rings, Beitr\"age Algebra Geom. {\bf 42} (2001), 395--400. 
      
\bibitem{Kadison}
L. Kadison, ``New examples of Frobenius extensions", {\sl University
Lect. Series} {\bf 14}, Amer. Math. Soc., Providence, 1999.

\bibitem{Nasta}
C. N\u{a}st\u{a}sescu and S. D\u{a}sc\u{a}lescu, Graded
$T$-rings, {\sl Comm. Algebra} {\bf 17} (1989), 3033--3042. 

\bibitem{Pareigis71}
B. Pareigis, When Hopf algebras are Frobenius Algebras, {\sl J.
Algebra} {\bf 18} (1971), 588--596.

\bibitem{Wisbauer}
R. Wisbauer, {On the category of comodules over corings},
in ``Mathematics and mathematics education (Bethlehem, 2000)", World Sci.
Publishing, River Edge, NJ, 2002, 325--336.


\end{thebibliography}
\end{document}